\documentclass[12pt]{article}

\usepackage[T1]{fontenc}
\usepackage{lmodern}
\usepackage{amsmath,amssymb,amsthm}
\usepackage{microtype}
\usepackage{needspace}
\usepackage[colorlinks=true,allcolors=blue,backref=page]{hyperref}

\makeatletter
\renewcommand{\section}{%
  \@startsection{section}{1}{\z@}%
    {-3.5ex \@plus -1ex \@minus -.2ex}%
    {2.3ex \@plus .2ex}%
    {\normalfont\fontfamily{cmr}\fontseries{bx}\Large\selectfont}}
\renewcommand{\subsection}{%
  \@startsection{subsection}{2}{\z@}%
    {-3.25ex \@plus -1ex \@minus -.2ex}%
    {1.5ex \@plus .2ex}%
    {\normalfont\fontfamily{cmr}\fontseries{bx}\large\selectfont}}
\makeatother

\theoremstyle{plain}
\newtheorem{theorem}{Theorem}[section]
\newtheorem{lemma}[theorem]{Lemma}

\newtheorem{corollary}[theorem]{Corollary}

\theoremstyle{definition}

\theoremstyle{remark}

\def\qed{\nopagebreak\hfill{\rule{4pt}{7pt}}
\medbreak}

\begin{document}
\begin{center}
  {\large\bfseries Congruences Modulo 16, 32 and 64 for a Restricted Overpartition Function \par}
  \vspace{12pt}
 Li Zhang\textsuperscript{1} and Qiao-Mei Zhang\textsuperscript{2}\\[6pt]
  \textsuperscript{1,2}Department of Mathematics and Physics,\\
  Shijiazhuang Tiedao University, Shijiazhuang 050043, Hebei, P. R. China \\[6pt]
  Email: \textsuperscript{1}\texttt{zhangli\_math7@stdu.edu.cn},
  \textsuperscript{2}\texttt{Zqm8970@163.com}
\end{center}

\noindent\textbf{Abstract.}
Let $\overline{t}^{(3)}(n)$ denote the number of overpartitions of $n$ in which
consecutive parts differ by a multiple of $3$ unless the larger part
is overlined, and the smallest part is overlined unless it is divisible by $3$.
Using dissections of theta functions, we establish two infinite families of congruences modulo
$32$ and $64$. We also obtain two residue formulas modulo $16$ and use them to derive four further infinite families of congruences modulo $16$. In particular, for every positive integer $M$ and nonnegative integer $n$, the residue of $\overline{t}^{(3)}(M^2(18n+12))$ modulo $16$ is $8$ if $6n+4$ is a square and $0$ otherwise.



\section{Introduction}
Recall that an overpartition of a nonnegative integer $n$ is a partition of $n$ where the last occurrence of each distinct part may be overlined.
In 2015, Bringmann et al. \cite{Bringmann-Dousse-Lovejoy-Mahlburg-2015} defined the partition function $\overline{t}^{(k)}(n)$ to count the number of overpartitions of $n$ satisfying (i) consecutive parts differ by a multiple of $k$ unless the larger of the two
is overlined, and (ii) the smallest part must be overlined unless it is divisible by $k$.

Throughout this paper, for each positive integer $k$, we write
\begin{equation*}
f_k:=(q^k;q^k)_{\infty}:=\prod_{n=1}^{\infty}(1-q^{kn}).
\end{equation*}
The combinatorial argument in Bringmann et al. \cite{Bringmann-Dousse-Lovejoy-Mahlburg-2015} yields
\begin{equation} \label{gf-tn}
\sum_{n=0}^{\infty}\overline{t}^{(k)}(n)q^n= \frac{f_{k+1}}{f_1f_k}.
\end{equation}

In recent years, the arithmetic properties of $\overline{t}^{(k)}(n)$ have attracted considerable attention. For $\overline{t}^{(2)}(n)$, extensive congruence results have been established; specifically, Hirschhorn and Sellers \cite{Hirschhorn-Sellers-2020}, and Chern and Hao \cite{Chern-Hao-2019} showed congruences modulo small powers of 2 and 3.
Lin et al. \cite{Lin-Liu-Wang-2020} proved two infinite families of congruences modulo 5 and 27.
The first author \cite{Zhang-2025, Zhang-2026} derived several infinite families of congruences modulo 81, 16, 32 and 64.
For other values of $k$, Naika and Gireesh \cite{Naika-Gireesh-2019} investigated  congruences modulo 2 and 4 for  $\overline{t}^{(3)}(n)$ and  $\overline{t}^{(7)}(n)$, modulo 4 and 5 for $\overline{t}^{(4)}(n)$, and modulo 3, 6 and 12 for  $\overline{t}^{(8)}(n)$.

In this paper, we use dissections of theta functions to study the arithmetic properties of $\overline{t}^{(3)}(n)$. Our first two results establish infinite families of congruences modulo $64$ and $32$ by extending the congruences on the progressions $144n+138$ and $144n+102$, respectively.

\begin{theorem}\label{thm-con-m2-144n+138}
Let $\alpha, n$ be nonnegative integers, and let $m$ be a positive integer with $\gcd(m,6)=1$. Then
\begin{align}\label{con-4a-m2-144n+138}
\overline{t}^{(3)}\big(4^{\alpha} m^2(144n+138)\big) & \equiv 0 \pmod{64},
\end{align}
and
\begin{align} \label{con-9a-m2-144n+138}
\overline{t}^{(3)}\big(9^{\alpha} m^2(144n+138)\big) & \equiv 0 \pmod{64}.
\end{align}
\end{theorem}

\begin{theorem}\label{thm-con-m2-144n+102}
Let $\alpha, n$ be nonnegative integers, and let $m$ be a positive integer with $\gcd(m,6)=1$. Then
\begin{align}
\overline{t}^{(3)}\big(4^{\alpha} m^2(144n+102)\big) & \equiv 0 \pmod{32}, \label{con-4a-m2-144n+102}
\end{align}
and
\begin{align}
\overline{t}^{(3)}\big(9^{\alpha} m^2(144n+102)\big) & \equiv 0 \pmod{32}. \label{con-9a-m2-144n+102}
\end{align}
\end{theorem}

We next determine the residues of $\overline{t}^{(3)}(n)$ modulo $16$ on square multiples of the progressions $18n+12$ and $24n+3$. These residue formulas yield four further infinite families of congruences indexed by primes.
For an odd prime $p$, $\big(\frac{\cdot}{p}\big)$ denotes the Legendre symbol.

\begin{theorem}\label{thm-rel-18n+12}
Let $\alpha, \beta, n$ be nonnegative integers, and let $m$ be a positive integer with $\gcd(m,6)=1$. Then
\begin{align}\label{rel-18n+12-end-1-mod16}
\overline{t}^{(3)}\big(4^{\alpha} 9^{\beta} m^2 (18n+12)\big)
\equiv
\begin{cases}
8, & \text{if $n = 6r^2+4r$ for some $r \in \mathbb{Z}$}, \\[5pt]
0, & \text{otherwise}
\end{cases} \pmod{16}.
\end{align}
Moreover, let $p\geq 5$ be a prime. For every integer $j$ with $1\leq j\leq p-1$, we have
\begin{align}\label{rel-18p^2+12-end-mod16}
\overline{t}^{(3)}\big(4^{\alpha} 9^{\beta} m^2 (18p^2n + 18pj + 12p^2)\big) \equiv 0 \pmod {16}.
\end{align}
For each such $j$ satisfying $\big(\frac{6j+4}{p} \big)=-1$, we also have
\begin{align}\label{rel-18p+12-end-mod16}
\overline{t}^{(3)}\big(4^{\alpha} 9^{\beta} m^2 (18pn + 18j + 12)\big) \equiv 0 \pmod {16}.
\end{align}
\end{theorem}

\begin{theorem}\label{thm-rel-24n+3}
Let $n$ be a nonnegative integer, and $m$ an odd integer. Then
\begin{align}\label{rel-m2-24n+3-mod16}
\overline{t}^{(3)}\big(m^2(24n+3)\big)  \equiv
\begin{cases}
4(-1)^{r+(m^2-1)/8}, & \text{if $n = 2r^2-r$ for some $r \in \mathbb{Z}$}, \\[5pt]
 0, & \text{otherwise}
\end{cases} \pmod{16}.
\end{align}
Moreover, let $p\geq 3$ be a prime. For every integer $j$ with $1\leq j\leq p-1$, we have
\begin{align}\label{rel-t24n+3-mod16}
\overline{t}^{(3)}\big(m^2(24p^2n+24pj+3p^2)\big) \equiv 0 \pmod {16}.
\end{align}
For each such $j$ satisfying $ \big( \frac{8j+1}{p} \big)=-1 $, we also have
\begin{align}\label{rel-24p+3-mod16}
\overline{t}^{(3)}\big(m^2(24pn+24j+3)\big) \equiv 0 \pmod {16}.
\end{align}
\end{theorem}

\section{Preliminaries}
In this section, we present some lemmas essential for proving the main theorems.
Recall that Ramanujan's theta functions $\psi(q)$ and $\Omega(q)$ are defined by
\begin{align}
\psi(q)
& := \sum_{n=-\infty}^{\infty} q^{2n^2-n}
   = \frac{f_2^2}{f_1},  \label{psi-def}  \\[5pt]
\Omega(q)
& := \sum_{n=-\infty}^{\infty}q^{3n^2+2n}
   = \frac{f_2^2 f_3 f_{12}}{f_1 f_4 f_6}. \label{Omega--def}
\end{align}
Replacing $q$ by $-q$  in \eqref{psi-def} gives
\begin{align}
\psi(-q)
& := \sum_{n=-\infty}^{\infty} (-1)^n q^{2n^2-n}
   = \frac{f_1f_4}{f_2}.  \label{psi-q-def}
\end{align}

\begin{lemma}\label{lem-f-2k}
For positive integers $k$ and $r$,
\begin{align*}
f_k^{2^r} \equiv f_{2k}^{2^{r-1}} \pmod {2^r}.
\end{align*}
\end{lemma}

\begin{lemma} The following 2-dissections hold.
\begin{align}
f_1^2 &= \frac{f_2f_8^5}{f_4^2f_{16}^2}-2q\frac{f_2f_{16}^2}{f_8},
\label{f1-2} \\[5pt]
\frac{1}{f_1^2} &= \frac{f_8^5}{f_2^5f_{16}^2}+2q\frac{f_4^2f_{16}^2}{f_2^5f_8}, \label{1-f1-2}  \\[5pt]
f_1^4 &= \frac{f_4^{10}}{f_2^2f_8^4} - 4q\frac{f_2^2f_8^4}{f_4^2}, \label{f1-4} \\[5pt]
\frac{1}{f_1^4} &= \frac{f_4^{14}}{f_2^{14}f_{8}^4}+4q\frac{f_4^2f_{8}^4}{f_2^{10}}, \label{1-f1-4} \\[5pt]
\frac{f_1^3}{f_3} &= \frac{f_4^3}{f_{12}} - 3q\frac{f_2^2 f_{12}^3}{f_4 f_6^2},
\label{f1-3-f3} \\[5pt]
\frac{f_1}{f_3^3}
&= \frac{f_2 f_4^2 f_{12}^2}{f_6^7} - q\frac{f_2^3 f_{12}^6}{f_4^2 f_6^9} ,
\label{f1-f3-3} \\[5pt]
\frac{f_3^3}{f_1}
&= \frac{f_4^3 f_6^2}{f_2^2 f_{12}} + q\frac{f_{12}^3}{f_4},
\label{f3-3-f1} \\[5pt]
\frac{f_3}{f_1} &= \frac{f_4 f_6 f_{16} f_{24}^2}{f_2^2 f_8 f_{12} f_{48}}
  + q\frac{f_6 f_8^2 f_{48}}{f_2^2 f_{16} f_{24}},   \label{f3-f1} \\[5pt]
f_1f_3 &= \frac{f_2f_8^2f_{12}^4}{f_4^2f_6f_{24}^2}
  - q\frac{f_4^4f_6f_{24}^2}{f_2f_8^2f_{12}^2}, \label{f1f3} \\[5pt]
\frac{1}{f_1f_3} &= \frac{f_8^2f_{12}^5}{f_2^2f_4f_6^4f_{24}^2}
  + q\frac{f_4^5f_{24}^2}{f_2^4f_6^2f_8^2f_{12}}. \label{1-f1f3}
\end{align}
\end{lemma}
{\noindent \textbf{Proof.}} Equations \eqref{f1-2}-\eqref{f3-f1} follows from \cite[Lemmas 2.2-2.6]{Yao-Xia-2013}. Equations \eqref{f1f3} and \eqref{1-f1f3} are given in \cite{Hirschhorn-2017} as Eqs. (30.12.1) and (30.12.3), respectively. \qed

\begin{lemma} The following 3-dissections hold.
\begin{align}
\frac{f_1^2}{f_2} &= \frac{f_9^2}{f_{18}} - 2q\frac{f_3f_{18}^2}{f_6f_9}, \label{f1-2-f2} \\[5pt]
\frac{f_2^2}{f_1} &= \frac{f_6f_9^2}{f_3f_{18}} + q\frac{f_{18}^2}{f_9}, \label{f2-2-f1} \\[5pt]
\frac{f_2}{f_1^2} &=\frac{f_6^4f_9^6}{f_3^8f_{18}^3} +2q\frac{f_6^3f_9^3}{f_3^7} +4q^2\frac{f_6^2f_{18}^3}{f_3^6}, \label{f2-f1-2} \\[5pt]
\frac{f_1}{f_2^2} &=\frac{f_3^2f_9^3}{f_6^6} - q\frac{f_3^3f_{18}^3}{f_6^7} + q^2\frac{f_3^4f_{18}^6}{f_6^8f_9^3}, \label{f1-f2-2} \\[5pt]
f_1^3 &= \frac{f_6f_9^6}{f_3f_{18}^3} - 3qf_9^3 + 4q^3\frac{f_3^2f_{18}^6}{f_6^2f_9^3},
\label{f1-3} \\[5pt]
\frac{f_4}{f_1} &= \frac{f_{12}f_{18}^4}{f_3^3f_{36}^2}
+ q\frac{f_6^2f_9^3f_{36}}{f_3^4f_{18}^2} + 2q^2\frac{f_6f_{18}f_{36}}{f_3^3},\label{f4-f1} \\[5pt]
f_1 f_2 &=
f_9 f_{18} \left( \frac{f_6 f_9^3}{f_3 f_{18}^3} - q - 2q^2 \frac{f_3 f_{18}^3}{f_6 f_9^3} \right).  \label{f1f2}
\end{align}
\end{lemma}

{\noindent \textbf{Proof.}}  Equations \eqref{f1-2-f2}-\eqref{f2-f1-2}, \eqref{f1-3} and \eqref{f1f2} follow directly from \cite[Eqs.(14.3.2), (14.3.3), (38.2.2), (14.8.5) and (14.3.1)]{Hirschhorn-2017}.
Equation \eqref{f1-f2-2} is established from \eqref{psi-def} and  \eqref{f2-2-f1}, while
\eqref{f4-f1} follows from \cite[Lemma 2.6, Eq. (57)]{Baruah-Ojah-2015}.
\qed

\begin{lemma}\label{lem-18n+12-tsim}   With the notation as above,
\begin{align*}
& \frac{f_2^3 f_3^2 f_6^4}{f_1^4f_4 f_{12}}
+ 2 \frac{f_3^4 f_6}{f_1^2}
- 2q \frac{f_1f_6^8}{f_2^3f_3}
- 2q^2 \frac{f_3^2f_6^{11}}{f_2^4}
\equiv \frac{f_3^4 f_6}{f_1^2}
- 2 \frac{f_4^4 f_6}{f_3^2}
+ 4 \frac{f_3^3 f_4^2}{f_1}  \pmod 8.
\end{align*}
\end{lemma}
{\noindent \textbf{Proof.}}  Set $X(q) = \frac{f_1^3}{f_3}$.  Using \eqref{f1-3-f3} and the identity $(-q;-q)_{\infty} = \frac{f_2^3}{f_1f_4}$, we obtain
\[
X(q) + X(-q) = \frac{f_1^3}{f_3} + \frac{f_2^9 f_3 f_{12}}{f_1^3f_4^3f_6^3} = 2\frac{f_4^3}{f_{12}}.
\]
Multiplying both sides by $\frac{f_3}{f_1^3 f_{12}}$, rearranging,  and applying Lemma \ref{lem-f-2k}  yields
\begin{align}\label{1-f12}
\frac{1}{f_{12}}
= - \frac{f_2^9 f_3^2}{f_1^6 f_4^3 f_6^3}
+ 2\frac{f_3 f_4^3}{f_1^3 f_{12}^2}
\equiv - \frac{f_1^2 f_3^2 f_4}{f_2^3 f_6^3}
+ 2\frac{f_1 f_3 f_4^3}{f_2^2 f_{6}^4} \pmod 8.
\end{align}
Substituting \eqref{1-f12} into the first term and  simplifying the resulting expression using  \eqref{f1-f3-3}, \eqref{f1f3} and Lemma \ref{lem-f-2k},  we derive
\begin{align*}
& \frac{f_2^3 f_3^2 f_6^4}{f_1^4f_4 f_{12}}
+ 2 \frac{f_3^4 f_6}{f_1^2}
- 2q \frac{f_1f_6^8}{f_2^3f_3}
- 2q^2 \frac{f_3^2f_6^{11}}{f_2^4}   \nonumber \\[5pt]
& ~~ \equiv  \frac{f_2^3 f_3^2 f_6^4}{f_1^4f_4} \left( \frac{1}{ f_{12}} \right)
+ 2 \frac{f_3^4 f_6}{f_1^2}
- 2q \frac{f_3^2 f_6^8}{f_2^3}
  \left(\frac{f_1}{f_3^3} + q \frac{f_6^{3}}{f_2} \right) \nonumber \\[5pt]
& ~~ \equiv \frac{f_3^4 f_6}{f_1^2}
+ 2 f_2^3 f_3^2 \left( f_1 f_3 - q \frac{f_6^5}{f_2} \right) \nonumber \\[5pt]
& ~~ \equiv \frac{f_3^4 f_6}{f_1^2}
- 2 \frac{f_4^4 f_6}{f_3^2}
+ 4 \frac{f_3^3 f_4^2}{f_1}  \pmod 8.
\end{align*}
 \qed

\section{Proof of Theorem \ref{thm-con-m2-144n+138}}
To prove Theorem \ref{thm-con-m2-144n+138}, we first establish the following specific congruence.
\begin{theorem}\label{thm-con-t144n+138}
For $n \geq 0$,
\begin{align}\label{con-t144n+138}
\overline{t}^{(3)}(144n+138) \equiv 0 \pmod {64}.
\end{align}
\end{theorem}
{\noindent \textbf{Proof.}}
From \eqref{gf-tn}, we get
\begin{align}\label{tn}
\sum_{n=0}^{\infty}\overline{t}^{(3)}(n)q^n
 = \frac{1}{f_3} \frac{f_4}{f_1}.
\end{align}
Plugging \eqref{f4-f1} into \eqref{tn}, extracting all terms of $q^{3n}$ and setting $q^3$ to $q$, we obtain
\begin{align}\label{t3n}
\sum_{n=0}^{\infty}\overline{t}^{(3)}(3n)q^n
= \frac{f_4 f_6^4}{f_{12}^2} \frac{1}{f_1^4}.
\end{align}
Substituting \eqref{1-f1-4} into \eqref{t3n}, extracting all terms of $q^{2n}$ and setting $q^2$ to $q$, we deduce
\begin{align}\label{t6n}
\sum_{n=0}^{\infty}\overline{t}^{(3)}(6n)q^n
& = \frac{f_2^{15}}{f_4^4 f_6^2} f_1^2 f_3^4 \Big( \frac{1}{f_1^4} \Big)^4.
\end{align}
Plugging  \eqref{f1-2}, \eqref{f1-4} and \eqref{1-f1-4} into \eqref{t6n}, extracting all terms of $q^{2n+1}$, dividing by $q$ and setting $q^2$ to $q$, and then applying Lemma \ref{lem-f-2k}, we derive
\begin{align}\label{t12n+6}
\sum_{n=0}^{\infty}\overline{t}^{(3)}(12n+6)q^n
& \equiv \left( 62 \frac{f_2^{52} f_6^{10} f_8^2}{f_1^{40} f_3^4 f_4^{17} f_{12}^4}
+ 16 \frac{f_2^{38} f_6^{10}}{f_1^{36} f_3^4 f_4^3 f_8^2 f_{12}^4} \right)
+ 60q \frac{f_2^{50} f_{12}^4}{f_1^{40} f_4^{11} f_6^2 f_8^2} \nonumber \\[5pt]
& \equiv 14 \frac{f_2^4 f_6^{10} f_8^2}{f_4 f_{12}^4} \frac{1}{f_3^4} \Big( \frac{1}{f_1^4} \Big)^2
+ 60q \frac{f_2^{10} f_{12}^4}{f_4^3 f_6^2 f_8^2} (f_1^4)^2   \pmod {64}.
\end{align}
Substituting \eqref{f1-4} and \eqref{1-f1-4} into \eqref{t12n+6}, extracting all terms of $q^{2n+1}$, dividing by $q$ and setting $q^2$ to $q$, and then applying Lemma \ref{lem-f-2k}, we find
\begin{align}\label{t24n+18}
\sum_{n=0}^{\infty}\overline{t}^{(3)}(24n+18)q^n
& \equiv 48 \frac{f_2^{15} f_4^2 f_6^{10}}{f_1^{20}f_3^4f_{12}^4}
 + 60 \frac{f_1^6 f_2^{17} f_6^4}{f_3^2 f_4^{10}}
 + 56q \frac{f_2^{27}f_{12}^4}{f_1^{24}f_4^6f_6^2}
   \nonumber \\[5pt]
& \equiv 48 f_2^5 f_4^2
+ 60 \frac{f_2 f_6^4}{f_4^2} \Big( \frac{f_1^3}{f_3} \Big)^2
+ 56q \frac{f_4^2 f_6^6}{f_2}   \pmod {64}.
\end{align}
Plugging  \eqref{f1-3-f3} into \eqref{t24n+18}, extracting all terms of $q^{2n+1}$,  dividing by $q$ and setting $q^2$ to $q$,  we are led to
\begin{align}\label{t48n+42}
\sum_{n=0}^{\infty}\overline{t}^{(3)}(48n+42)q^n
\equiv 24 f_3^2 f_6^2 f_1^3
+ 56 f_3^6 \frac{f_2^2}{f_1}  \pmod {64}.
\end{align}
Substituting \eqref{f2-2-f1}  and \eqref{f1-3} into \eqref{t48n+42}, we arrive at
\begin{align}\label{t48n+42-1}
 \sum_{n=0}^{\infty}\overline{t}^{(3)}(48n+42)q^n
& \equiv 24 \frac{f_3 f_6^3 f_9^6}{f_{18}^3}
+ 56 \frac{f_3^5 f_6 f_9^2}{f_{18}}
+ 56q \frac{f_3^6 f_{18}^2}{f_9} \nonumber \\[5pt]
&~~ + 56q f_3^2 f_6^2 f_9^3
+ 32q^3 \frac{f_3^4 f_{18}^6}{f_9^3}  \pmod {64}.
\end{align}
Since the right-hand side of \eqref{t48n+42-1} has no terms of the form $q^{3n+2}$, it implies \eqref{con-t144n+138}.  \qed

{\noindent \textbf{Proof of Theorem \ref{thm-con-m2-144n+138}.}}
Rewriting \eqref{t6n} as
\begin{align}\label{t6n-1}
\sum_{n=0}^{\infty}\overline{t}^{(3)}(6n)q^n
& = \frac{f_3^4}{f_6^2} \Big( \frac{f_2}{f_1^2} \Big)^7 \Big( \frac{f_2^2}{f_4} \Big)^4.
\end{align}
Plugging  \eqref{f1-2-f2} and \eqref{f2-f1-2} into \eqref{t6n-1}, extracting all terms of $q^{3n+2}$, dividing by $q^2$ and setting $q^3$ to $q$, and then applying Lemma \ref{lem-f-2k} and Lemma \ref{lem-18n+12-tsim}, we obtain
\begin{align}\label{t18n+12}
\sum_{n=0}^{\infty}\overline{t}^{(3)}(18n+12)q^n
& \equiv  56\frac{f_2^{27} f_3^{42}}{f_1^{52} f_4 f_6^{16} f_{12}}
+ 48\frac{f_2^{24} f_3^{36}}{f_1^{50} f_6^7 f_{12}^4}
+ 16q \frac{f_2^{27} f_3^{39} f_{12}^2}{f_1^{51} f_4^2 f_6^{16}}
+ 16q^2 \frac{f_2^{30} f_3^{42} f_{12}^8}{f_1^{52} f_4^4 f_6^{25}} \nonumber \\[5pt]
& \equiv  -8 \left\{ \frac{f_2^3 f_3^2 f_6^4}{f_1^4f_4 f_{12}}
+ 2 \frac{f_3^4 f_6}{f_1^2}
- 2q \frac{f_1f_6^8}{f_2^3f_3}
- 2q^2 \frac{f_3^2f_6^{11}}{f_2^4}  \right \} \nonumber \\[5pt]
& \equiv -8 f_6 f_3^4 \frac{1}{f_1^2}
+ 16 f_4^4 f_6 \frac{1}{f_3^2}
+ 32 f_4^2 \frac{f_3^3}{f_1}     \pmod {64}.
\end{align}
Substituting \eqref{1-f1-2}, \eqref{f1-4} and \eqref{f3-3-f1} into \eqref{t18n+12}, extracting all terms of $q^{2n}$ and setting $q^2$ to $q$, and then applying Lemma \ref{lem-f-2k}, we derive
\begin{align}\label{t36n+12}
\sum_{n=0}^{\infty}\overline{t}^{(3)}(36n+12)q^n
& \equiv -8 \frac{f_4^5 f_6^{10}}{f_1^5 f_3 f_8^2 f_{12}^4}
+ 16 \frac{f_2^4 f_{12}^5}{f_3^4 f_{24}^2}
+ 32 \frac{f_2^5 f_3^2}{f_1^2 f_6}
\nonumber \\[5pt]
&\equiv - 8\frac{f_4^5 f_6^2}{f_8^2} \frac{1}{f_1^4} \frac{1}{f_1 f_3}
+ 16 \frac{f_4^2 f_{12}}{f_6^2} + 32 f_4^2  \pmod {64}.
\end{align}
Plugging \eqref{1-f1-4} and \eqref{1-f1f3} into \eqref{t36n+12}, extracting all terms of $q^{2n+1}$, dividing by $q$ and setting $q^2$ to $q$, and then applying Lemma \ref{lem-f-2k}, we get
\begin{align}\label{t72n+48}
\sum_{n=0}^{\infty}\overline{t}^{(3)}(72n+48)q^n
& \equiv -8 \frac{f_2^{24}f_{12}^2}{f_1^{18} f_4^8 f_6}
 + 32 \frac{f_2^6 f_4^4 f_6^5}{f_1^{12} f_3^2 f_{12}^2} \nonumber \\[5pt]
&\equiv -8 \frac{f_{12}^2}{f_1^2 f_6} + 32 f_4^4   \pmod {64}.
\end{align}
Subtracting \eqref{t72n+48} from \eqref{t18n+12}, applying \eqref{1-f1-4} and Lemma \ref{lem-f-2k}, employing \eqref{1-f1-2}, \eqref{f3-3-f1} and \eqref{f3-f1},  and  then invoking Lemma \ref{lem-f-2k} again, we find
\begin{align}\label{t18n+12-t72n+48}
& \sum_{n=0}^{\infty}\Big(\overline{t}^{(3)}(18n+12)
- \overline{t}^{(3)}\big(4(18n+12)\big)\Big)q^n \nonumber \\[5pt]
&~~~ \equiv 8 \frac{f_3^4 f_{12}^2}{f_1^2 f_6} \left( \frac{1}{f_3^4}
- \frac{f_6^2}{f_{12}^2} \right)
+ \left( 16 \frac{f_4^4 f_6}{f_3^2} - 32 f_4^4 \right)
+ 32\frac{f_3^3 f_4^2}{f_1} \nonumber \\[5pt]
&~~~ \equiv -16 f_8^2 f_6 \frac{1}{f_3^2}
+ 32f_8 \frac{f_3^3}{f_1}
+ 32q^3 f_{12}^7 \frac{f_6}{f_2}   \nonumber \\[5pt]
&~~~ \equiv - 16 \frac{f_8^2 f_{24}}{f_{12}^2}
+ 32 f_8^2 + 32 q \frac{f_8 f_{12}^3}{f_4}
+ 32 q^5 \frac{f_{24}^6}{f_8}   \pmod {64}.
\end{align}
Since the right-hand side of \eqref{t18n+12-t72n+48} has no terms of the form $q^{4n+2}$,  it follows
\begin{align}\label{t4(18n+12)-t4^2(18n+12)}
\overline{t}^{(3)}\big(4^2(18n+12)\big) \equiv \overline{t}^{(3)}\big(4(18n+12)\big)  \pmod {64}.
\end{align}
For $\alpha \geq 0$, setting $n \rightarrow 4^{\alpha}\cdot 8n + \frac{4^{\alpha} \cdot 23 - 2}{3}$ in \eqref{t4(18n+12)-t4^2(18n+12)}, we arrive at
\begin{align}\label{rel-4a-144n+138-rel}
\overline{t}^{(3)}\big(4^{\alpha+2}(144n+138)\big)
\equiv \overline{t}^{(3)}\big(4^{\alpha + 1}(144n+138)\big)   \pmod{64}.
\end{align}
Also \eqref{t18n+12-t72n+48} contains no terms of $q^{8n+7}$, which leads to
\begin{align}\label{rel-4-144n+138}
\overline{t}^{(3)}\big(4(144n+138)\big)
\equiv \overline{t}^{(3)}(144n+138)   \pmod{64}.
\end{align}
Combining \eqref{rel-4a-144n+138-rel}, \eqref{rel-4-144n+138}  and \eqref{con-t144n+138}, we conclude that for $\alpha \geq 0$,
\begin{align}\label{con-4a-t144n+138}
\overline{t}^{(3)}\big(4^{\alpha}(144n+138)\big) \equiv 0 \pmod {64}.
\end{align}

On the other hand, extracting all terms of $q^{3n+1}$ in \eqref{t48n+42-1}, dividing by $q$ and setting $q^3$ to $q$,  we derive
\begin{align}\label{t144n+90}
\sum_{n=0}^{\infty}\overline{t}^{(3)}(144n+90)q^n
\equiv 56 \frac{f_6^2}{f_3} (f_1^3)^2
+ 56 f_3^3 (f_1 f_2)^2  \pmod {64}.
\end{align}
Plugging \eqref{f1-3} and \eqref{f1f2} into \eqref{t144n+90}, extracting all terms of $q^{3n+2}$, dividing by $q^2$ and setting $q^3$ to $q$,  we deduce
\begin{align}\label{t9(48n+42)}
\sum_{n=0}^{\infty}\overline{t}^{(3)}\big(9(48n+42)\big)q^n
& \equiv 24 f_3^2 f_6^2 f_1^3
+ 56 f_3^6 \frac{f_2^2}{f_1}   \pmod {64}.
\end{align}
In view of \eqref{t48n+42} and \eqref{t9(48n+42)}, we get
\begin{align}\label{t48n+42-t9(48n+42)}
\overline{t}^{(3)}(48n+42) \equiv \overline{t}^{(3)}\big(9(48n+42)\big)  \pmod {64}.
\end{align}
For $\alpha \geq 0$, setting
$n\rightarrow 9^{\alpha} \cdot 3n + \frac{9^{\alpha}\cdot23-7}{8}$
in \eqref{t48n+42-t9(48n+42)}, we obtain
\begin{align}\label{t9a(48n+42)-t9a(144n+138)}
\overline{t}^{(3)}\big(9^{\alpha}(144n+138)\big) \equiv \overline{t}^{(3)}\big(9^{\alpha+1}(144n+138)\big)  \pmod {64}.
\end{align}
Combining \eqref{con-t144n+138}  and \eqref{t9a(48n+42)-t9a(144n+138)}, we arrive at for $\alpha \geq 0$,
\begin{align}\label{con-9a-144n+138}
\overline{t}^{(3)}\big(9^{\alpha}(144n+138)\big) \equiv 0 \pmod {64}.
\end{align}
Since $\gcd(m,6)=1$, we have $m^2\equiv1\pmod{24}$. For $c\in \{17, 23\}$, define
\[
n_c=m^2n+\frac{c(m^2-1)}{24}.
\]
Then $n_c$ is a nonnegative integer and $m^2(144n+6c)=144n_c+6c$. Replacing $n$ by $n_{23}$ in (3.17) and (3.22) gives (1.2) and (1.3), respectively. This completes the proof.
\qed

\section{Proof of Theorem \ref{thm-con-m2-144n+102}}
{\noindent \textbf{Proof of Theorem \ref{thm-con-m2-144n+102}.}}
Reducing \eqref{t18n+12} modulo 32 and applying Lemma \ref{lem-f-2k}, then employing \eqref{1-f1-2}, \eqref{f3-3-f1} and reapplying Lemma \ref{lem-f-2k}, we obtain
\begin{align}\label{t18n+12-mod32}
\sum_{n=0}^{\infty}\overline{t}^{(3)}(18n+12)q^n
& \equiv - 8 \frac{f_3^4 f_6}{f_1^2} + 16\frac{f_4^4 f_6}{f_3^2}
 \equiv - 8 f_2 \frac{1}{f_1^2}\frac{f_6^3}{f_2} + 16 f_4^4   \nonumber \\[5pt]
& \equiv  24 \frac{f_8^2 f_{12}^2}{f_{24}} + 16 f_4^4
+ 16 q f_8^5 + 24 q^2 \frac{f_{24}^3}{f_4^2}
+ 16 q^3 f_{16} f_{24}^3   \pmod {32}.
\end{align}
Note that there are no terms of $q^{8n+5}$ in \eqref{t18n+12-mod32}. Hence it follows
\begin{align}\label{con-144n+102-mod32}
\overline{t}^{(3)}(144n+102)\equiv 0  \pmod {32}.
\end{align}
Taking  \eqref{t18n+12-t72n+48} modulo 32 and applying Lemma \ref{lem-f-2k}, we find
\begin{align}\label{t18n+12-t72n+48-mod32}
\sum_{n=0}^{\infty}\Big(\overline{t}^{(3)}(18n+12)
- \overline{t}^{(3)}\big(4(18n+12)\big)\Big)q^n
 \equiv  16 \frac{f_8^2 f_{24}}{f_{12}^2}
 \equiv  16 f_8^2   \pmod {32}.
\end{align}
Clearly, \eqref{t18n+12-t72n+48-mod32} contains no terms of the form $q^{8n+5}$, which implies
\begin{align}\label{rel-144n+102-4(144n+102)-mod32}
\overline{t}^{(3)}(144n+102)
 \equiv  \overline{t}^{(3)}\big(4(144n+102)\big)   \pmod {32}.
\end{align}
Reducing \eqref{t4(18n+12)-t4^2(18n+12)} modulo $32$ and replacing $n$ by $4^{\alpha}\cdot 8n + \frac{4^{\alpha} \cdot 17 - 2}{3}$, we obtain, for every nonnegative integer $\alpha \geq 0$,
\begin{align}\label{rel-4a-144n+138-mod32}
\overline{t}^{(3)}\big(4^{\alpha+2}(144n+102)\big)
\equiv \overline{t}^{(3)}\big(4^{\alpha + 1}(144n+102)\big)   \pmod{32}.
\end{align}
Combining \eqref{con-144n+102-mod32}, \eqref{rel-144n+102-4(144n+102)-mod32} and \eqref{rel-4a-144n+138-mod32}, we arrive at for $\alpha \geq 0$,
\begin{align}\label{t4a(144n+102)-mod32}
\overline{t}^{(3)}\big(4^{\alpha}(144n+102)\big) \equiv 0  \pmod {32}.
\end{align}

On the other hand, invoking Lemma \ref{lem-f-2k}, we rewrite \eqref{t6n} as
\begin{align}\label{t6n-mod32}
\sum_{n=0}^{\infty}\overline{t}^{(3)}(6n)q^n
\equiv  \frac{f_1^{18} f_3^4}{f_2 f_4^4 f_6^2}
\equiv \frac{f_3^4}{f_6^2} \Big( \frac{f_1^2}{f_2} \Big)^9 \Big( \frac{f_2^2}{f_4} \Big)^4  \pmod {32}.
\end{align}
Substituting \eqref{f1-2-f2} into \eqref{t6n-mod32}, extracting all terms of $q^{3n}$ and setting $q^3$ to $q$, and then applying Lemma \ref{lem-f-2k}, we obtain
\begin{align}\label{t18n-mod32}
\sum_{n=0}^{\infty}\overline{t}^{(3)}(18n)q^n
&\equiv \frac{f_1^4 f_3^{18}}{f_2^2 f_6 f_{12}^4}
+ 16q \frac{f_1^5 f_3^{15}}{f_2^2 f_4 f_6 f_{12}}  \nonumber \\[5pt]
&\equiv \frac{f_3^{18}}{f_6 f_{12}^4} \Big( \frac{f_1^2}{f_2} \Big)^2
+ 16q f_3^9 \frac{f_1}{f_2^2}   \pmod {32}.
\end{align}
Plugging \eqref{f1-2-f2} and \eqref{f1-f2-2} into \eqref{t18n-mod32}, extracting all terms of $q^{3n}$ and setting $q^3$ to $q$, and then applying Lemma \ref{lem-f-2k}, we deduce
\begin{align}\label{t54n-mod32}
\sum_{n=0}^{\infty}\overline{t}^{(3)}(54n)q^n
&\equiv  \frac{f_1^{18} f_3^4}{f_2 f_4^4 f_6^2} + 16q \frac{f_1^{13}f_6^6}{f_2^8 f_3^3} \nonumber \\[5pt]
&\equiv \frac{f_3^4}{f_6^2} \Big( \frac{f_1^2}{f_2} \Big)^9 \Big( \frac{f_2^2}{f_4} \Big)^4
+ 16q f_3^9 \frac{f_1}{f_2^2}   \pmod {32}.
\end{align}
Substituting \eqref{f1-2-f2} and \eqref{f1-f2-2} into \eqref{t54n-mod32}, extracting all terms of $q^{3n}$ and setting $q^3$ to $q$, and then applying Lemma \ref{lem-f-2k}, we derive
\begin{align}\label{t162n-mod32}
\sum_{n=0}^{\infty}\overline{t}^{(3)}(162n)q^n
& \equiv \frac{f_1^4 f_3^{18}}{f_2^2 f_6 f_{12}^4}
+ \left( 16q \frac{f_1^{13} f_6^6}{f_2^8 f_3^3}
+ 16q \frac{f_1^5 f_3^{15}}{f_2^2 f_4 f_6 f_{12}} \right) \nonumber \\[5pt]
& \equiv \frac{f_3^{18}}{f_6 f_{12}^4} \Big( \frac{f_1^2}{f_2} \Big)^2   \pmod {32}.
\end{align}
Plugging \eqref{f1-2-f2} into \eqref{t162n-mod32}, extracting all terms of $q^{3n}$ and setting $q^3$ to $q$, we find
\begin{align}\label{t81*6n-mod32}
\sum_{n=0}^{\infty}\overline{t}^{(3)}(81 \cdot 6n)q^n
\equiv \frac{f_1^{18} f_3^4}{f_2 f_4^4 f_6^2}    \pmod {32}.
\end{align}
Combining \eqref{t6n-mod32} and \eqref{t81*6n-mod32},  we are led to \begin{align*}
\overline{t}^{(3)}(6n) \equiv \overline{t}^{(3)}(81 \cdot 6n)  \pmod {32}.
\end{align*}
It follows that for $\alpha \geq 0$,
\begin{align}\label{rel-6n-9^2a*6n}
\overline{t}^{(3)}(6n) \equiv \overline{t}^{(3)}(9^{2\alpha} \cdot 6n)  \pmod {32}.
\end{align}
Substituting \eqref{f1-2-f2} and \eqref{f1-f2-2} into \eqref{t54n-mod32}, extracting all terms of $q^{3n+2}$, dividing by $q^2$ and setting $q^3$ to $q$, and then applying Lemma \ref{lem-f-2k} and Lemma  \ref{lem-18n+12-tsim}, we derive
\begin{align}\label{t162n+108-mod32}
& \sum_{n=0}^{\infty}\overline{t}^{(3)}(162n+108)q^n  \nonumber \\[5pt]
&~~~ \equiv \left( - 8 \frac{f_1^4 f_3^{18}}{f_2 f_4 f_6^4 f_{12}}
+ 16\frac{f_1^6 f_3^{12} f_6^5}{f_2^4 f_{12}^4}
+ 16q \frac{f_1^5 f_3^{15} f_{12}^2}{f_2 f_4^2 f_6^4}
+ 16q^2 \frac{f_1^4 f_2^2 f_3^{18} f_{12}^8}{f_4^4 f_6^{13}}  \right)
+ 16 \frac{f_1^{12} f_6^3}{f_2^7}   \nonumber \\[5pt]
&~~~ \equiv  -8\left( \frac{f_2^3 f_3^2 f_6^4}{f_1^4f_4 f_{12}}
+ 2 \frac{f_3^4 f_6}{f_1^2}
- 2q \frac{f_1f_6^8}{f_2^3f_3}
- 2q^2 \frac{f_3^2f_6^{11}}{f_2^4} \right)
+ 16 \frac{f_6^3}{f_1^2}   \nonumber \\[5pt]
&~~~ \equiv 8 f_6^3 \frac{1}{f_1^2} + 16 f_4^4    \pmod {32}.
\end{align}
Plugging \eqref{1-f1-2} into \eqref{t162n+108-mod32}, extracting all terms of $q^{2n+1}$, dividing by $q$ and setting $q^2$ to $q$, applying Lemma \ref{lem-f-2k}, employing \eqref{f3-3-f1} and reapplying Lemma \ref{lem-f-2k},  we get
\begin{align}\label{t9(36n+30)-mod32}
\sum_{n=0}^{\infty}\overline{t}^{(3)}\big(9(36n+30)\big)q^n
& \equiv 16 \frac{f_2^2 f_3^3 f_8^2}{f_1^5 f_4}
 \equiv 16 f_4^3 \frac{f_3^3}{f_1}  \nonumber \\[5pt]
&\equiv 16 f_4^5 + 16 q f_8 f_{12}^3
 \pmod {32}.
\end{align}
In \eqref{t18n+12-mod32}, extracting all terms of $q^{2n+1}$, dividing by $q$ and setting $q^2$ to $q$, we obtain
\begin{align}\label{t(36n+30)-mod32}
\sum_{n=0}^{\infty}\overline{t}^{(3)}(36n+30)q^n
&\equiv 16 f_4^5 + 16 q f_8 f_{12}^3   \pmod {32}.
\end{align}
Comparing \eqref{t9(36n+30)-mod32} and \eqref{t(36n+30)-mod32} gives
\begin{align}\label{t36n+30-t9(36n+30)-mod32}
\overline{t}^{(3)}(36n+30) \equiv \overline{t}^{(3)}\big(9(36n+30)\big)   \pmod {32}.
\end{align}
Substituting $n\rightarrow 6n+5$ and $n\rightarrow 9(6n+5)$ into \eqref{rel-6n-9^2a*6n}, and  using \eqref{t36n+30-t9(36n+30)-mod32}, we conclude
\begin{align}\label{t36n+30-t9a(36n+30)-mod32}
\overline{t}^{(3)}(36n+30) \equiv \overline{t}^{(3)}\big(9^{\alpha}(36n+30)\big)   \pmod {32}.
\end{align}
Setting $n\rightarrow 4n+2$ in \eqref{t36n+30-t9a(36n+30)-mod32} and then combining with \eqref{con-144n+102-mod32} yields
\begin{align}\label{t9a(144n+102)-mod32}
\overline{t}^{(3)}\big(9^{\alpha}(144n+102)\big) \equiv 0  \pmod {32}.
\end{align}
Using the substitution from Section 3 with $c=17$, we have $m^2(144n+102)=144n_{17}+102$. Replacing $n$ by $n_{17}$ in \eqref{t4a(144n+102)-mod32} and \eqref{t9a(144n+102)-mod32} gives \eqref{con-4a-m2-144n+102} and \eqref{con-9a-m2-144n+102}, respectively.  \qed

\section{Proof of Theorem \ref{thm-rel-18n+12}}
{\noindent \textbf{Proof of Theorem \ref{thm-rel-18n+12}.}}
Reducing \eqref{t18n+12} modulo 16, then applying Lemma \ref{lem-f-2k} and \eqref{Omega--def}, we derive
\begin{align}\label{t18n+12-mod16}
\sum_{n=0}^{\infty}\overline{t}^{(3)}(18n+12)q^n
\equiv 8\frac{f_6^3}{f_2}
\equiv 8 \Omega(q^2)
= 8\sum_{r=-\infty}^{\infty} q^{6r^2+4r} \pmod{16}.
\end{align}
Reducing \eqref{t18n+12-t72n+48} modulo $16$  gives
\begin{align}\label{rel-t18n+12-t72n+48-mod16}
\overline{t}^{(3)}(18n+12) \equiv \overline{t}^{(3)}\big(4(18n+12)\big) \pmod {16}.
\end{align}
Since $4(18n+12)=18(4n+2)+12$, iterating \eqref{rel-t18n+12-t72n+48-mod16}  yields
\begin{align}\label{rel-t4a(18n+12)-mod16}
\overline{t}^{(3)}(18n+12)
\equiv \overline{t}^{(3)}\big(4^{\alpha} (18n+12)\big) \pmod {16},
\end{align}
for every nonnegative integer $\alpha$.
By \eqref{t6n-mod32} and \eqref{t54n-mod32}, we find
\begin{align}\label{rel-t6n-t54n-mod16}
\overline{t}^{(3)}(6n) \equiv \overline{t}^{(3)}(54n)    \pmod {16}.
\end{align}
Iterating \eqref{rel-t6n-t54n-mod16} produces
\begin{align}\label{rel-t6n-t9a-6n-mod16}
\overline{t}^{(3)}(6n) \equiv \overline{t}^{(3)}(9^{\beta} \cdot 6n) \pmod {16},
\end{align}
for every nonnegative integer $\beta$.
Setting $n\rightarrow 4^{\alpha}(3n+2)$ in \eqref{rel-t6n-t9a-6n-mod16} and then combining with \eqref{rel-t4a(18n+12)-mod16},  we obtain
\begin{align}\label{t4a-9b(18n+12)-mod16}
\overline{t}^{(3)}(18n+12) \equiv \overline{t}^{(3)}\big(4^{\alpha} 9^{\beta}(18n+12)\big)  \pmod {16}.
\end{align}
Combining \eqref{t18n+12-mod16} and \eqref{t4a-9b(18n+12)-mod16}, we arrive at
\begin{align}\label{t4a-9b(18n+12)-ar-mod16}
\overline{t}^{(3)}\big(4^{\alpha} 9^{\beta}(18n+12)\big)
\equiv  \left\{
          \begin{array}{ll}
            8 \pmod {16}, & \hbox{if $n=6r^2+4r$ for some integer $r$;} \\[5pt]
            0 \pmod {16}, & \hbox{otherwise.}
          \end{array}
        \right.
\end{align}
Since $\gcd(m,6)=1$,  $\frac{2(m^2-1)}{3}$ is an integer. And
\begin{align}\label{t4a-9b(18n+12)-m2-mod16}
4^{\alpha} 9^{\beta}m^2(18n+12)
= 4^{\alpha} 9^{\beta}\left(18\Big(m^2n + \frac{2(m^2-1)}{3}\Big)+12\right).
\end{align}
For some integer $r$,
\begin{align}\label{n=2r(3r+2)}
n = 6r^2+4r
& \Leftrightarrow 6n+4 = (6r+2)^2 \nonumber \\[5pt]
& \Leftrightarrow m^2(6n+4) = m^2(6r+2)^2 \nonumber \\[5pt]
& \Leftrightarrow m^2n+ \frac{2(m^2-1)}{3} = 6s^2+4s,
\end{align}
where we set $s:= \frac{m-1}{3} + mr$ if $m\equiv 1 \pmod 6$, and $s:= \frac{-m-1}{3} - mr$ if $m\equiv -1 \pmod 6$.
Combining \eqref{t4a-9b(18n+12)-ar-mod16}-\eqref{n=2r(3r+2)}, we deduce \eqref{rel-18n+12-end-1-mod16}.

Let $p$ be an odd prime and let $a, b, n, j$ be integers with $p\nmid a$. If $p\nmid j$, then
\[
ap^2n + apj + bp^2 = p\big(p(an+b)+aj\big)
\]
has $p$-adic valuation $1$ and hence cannot be a square. If $\big(\frac{aj+b}{p}\big)=-1$, then $a(pn+j)+b$ is a quadratic nonresidue modulo $p$ and hence cannot be a square.

To prove \eqref{rel-18p^2+12-end-mod16} and  \eqref{rel-18p+12-end-mod16}, put
\[
u = p^2n+pj+\frac{2(p^2-1)}{3},  \qquad   v=pn+j.
\]
Both indices are nonnegative integers under the stated hypotheses, since \(p\ge5\) is prime and $p^2\equiv1\pmod3$. With $(a,b)=(6,4)$, the two square tests above apply to
\[
6u+4=6p^2n+6pj+4p^2,  \qquad  6v+4=6(pn+j)+4,
\]
under the hypotheses of \eqref{rel-18p^2+12-end-mod16} and  \eqref{rel-18p+12-end-mod16}, respectively. Applying the zero case of \eqref{rel-18n+12-end-1-mod16} with $n$ replaced by $u$ and $v$ gives the two congruences.
\qed

\begin{corollary}
Let $\alpha, \beta, n$ be nonnegative integers, and let $m$ be a positive integer with $\gcd(m,6)=1$. Then
\begin{align}
\overline{t}^{(3)}\big(4^{\alpha} 9^{\beta} m^2 (36n+30)\big)
& \equiv 0 \pmod {16}, \label{t4a-9b-m2(36n+30)-mod16} \\[5pt]
\overline{t}^{(3)}\big(4^{\alpha} 9^{\beta} m^2 (144n+84)\big)
& \equiv 0 \pmod {16}, \label{t4a-9b-m2(144n+84)-mod16} \\[5pt]
\overline{t}^{(3)}\big(4^{\alpha} 9^{\beta} m^2 (144n+120)\big)
& \equiv 0 \pmod {16}, \label{t4a-9b-m2(144n+120)-mod16} \\[5pt]
\overline{t}^{(3)}\big(4^{\alpha} 9^{\beta} m^2 (288n+156)\big)
& \equiv 0 \pmod {16}. \label{t4a-9b-m2(288n+156)-mod16}
\end{align}
\end{corollary}
{\noindent \textbf{Proof.}}
Since $6r^2+4r \equiv 0, 2, 10 \pmod{16}$, so the substitutions $n\rightarrow 2n+1, 8n+4, 8n+6, 16n+8$ all fall into the zero case of (1.6) and yield
\eqref{t4a-9b-m2(36n+30)-mod16}-\eqref{t4a-9b-m2(288n+156)-mod16}, respectively.
\qed

\section{Proof of Theorem \ref{thm-rel-24n+3}}
{\noindent \textbf{Proof of Theorem \ref{thm-rel-24n+3}.}}
Plugging \eqref{1-f1-4} into \eqref{t3n}, extracting all terms of $q^{2n+1}$, dividing by $q$ and setting $q^2$ to $q$, and then applying Lemma \ref{lem-f-2k}, we deduce
\begin{align}\label{t6n+3-mod16}
\sum_{n=0}^{\infty}\overline{t}^{(3)}(6n+3)q^n
= 4\frac{f_2^3 f_3^4 f_4^4}{f_1^{10} f_6^2}
\equiv 4\frac{f_4^4}{f_2} \frac{1}{f_1^2} \pmod {16}.
\end{align}
Substituting \eqref{1-f1-2} into \eqref{t6n+3-mod16}, extracting all terms of $q^{2n}$ and setting $q^2$ to $q$, and then applying Lemma \ref{lem-f-2k}, we derive
\begin{align}\label{t12n+3-mod16}
\sum_{n=0}^{\infty}\overline{t}^{(3)}(12n+3)q^n
\equiv 4\frac{f_2^4 f_4^5}{f_1^6 f_8^2}
\equiv 4f_4 f_1^2  \pmod {16}.
\end{align}
Plugging \eqref{f1-2} into \eqref{t12n+3-mod16}, extracting all terms of $q^{2n}$ and setting $q^2$ to $q$, and then applying Lemma \ref{lem-f-2k} and \eqref{psi-q-def}, we obtain
\begin{align*}
\sum_{n=0}^{\infty}\overline{t}^{(3)}(24n+3)q^n
\equiv 4\frac{f_1 f_4^5}{f_2 f_8^2}
\equiv 4\psi(-q)
= 4\sum_{r=-\infty}^{\infty} (-1)^r q^{2r^2-r}   \pmod {16}.
\end{align*}
It follows that
\begin{align}\label{rel-24n+3}
\overline{t}^{(3)}(24n+3)
\equiv
\begin{cases}
4(-1)^r \pmod{16}, & \text{if $n = 2r^2-r$ for some integer $r$};\\[5pt]
~~~~~~ 0 \ \pmod{16}, & \text{otherwise}.
\end{cases}
\end{align}
Since $m$ is an odd integer, $\frac{m^2-1}{8}$ is an integer. And
\begin{align}\label{rel-m2-24n+3}
m^2(24n+3)=24 \left( m^2n+ \frac{m^2-1}{8} \right) + 3.
\end{align}
For some integer $r$,
\begin{align}\label{n=r(2r-1)}
n = 2r^2-r
& \Leftrightarrow 8n+1 = (4r-1)^2 \nonumber \\[5pt]
& \Leftrightarrow m^2(8n+1) = m^2(4r-1)^2 \nonumber \\[5pt]
& \Leftrightarrow m^2n+ \frac{m^2-1}{8} = 2s^2-s,
\end{align}
where we set $s:= \frac{-m+1}{4} + mr$ if $m\equiv 1 \pmod 4$, and $s:= \frac{m+1}{4} - mr$ if $m\equiv -1 \pmod 4$.
Combining \eqref{rel-24n+3}-\eqref{n=r(2r-1)}, we deduce \eqref{rel-m2-24n+3-mod16}.

To prove \eqref{rel-t24n+3-mod16} and \eqref{rel-24p+3-mod16},  put
\[
u=p^2n+pj+\frac{p^2-1}{8},
\qquad v=pn+j.
\]
Both indices are nonnegative integers under the stated hypotheses, since $p$ is odd and $p^2\equiv1\pmod8$. Applying the square tests from Section 5 with $(a,b)=(8,1)$, we find that
\[
8u+1 = 8p^2n+8pj+p^2,
\qquad
8v+1 = 8(pn+j)+1
\]
cannot be squares under the hypotheses of \eqref{rel-t24n+3-mod16} and \eqref{rel-24p+3-mod16}, respectively. Applying the zero case of \eqref{rel-m2-24n+3-mod16} with $n$ replaced by $u$ and $v$ gives the two congruences.
\qed

\vspace{.3cm}

{\noindent \bf Acknowledgments}
This work was supported by the National Natural Science Foundation of China (Grant No. 11901403) and the Fundamental Research Funds for Shijiazhuang Tiedao University (Grant No. ZQK202610).

\vspace{.1cm}
{\noindent \bf Data Availability}
All data generated  during this study are included in this article.

\vspace{.1cm}
{\noindent \bf Conflicts of Interest}
The authors declare no conflict of interest.

\end{document}